\documentclass{article}
\usepackage{amsmath}
\usepackage[dvips]{graphicx}
\usepackage{psfrag}
\newcommand{\openbox}{\leavevmode
  \hbox to.77778em{%
  \hfil\vrule
  \vbox to.675em{\hrule width.6em\vfil\hrule}%
  \vrule\hfil}}
\newenvironment{proof}[1][Proof]
{\par\addvspace{6pt}\normalfont \itshape #1\@{.}\hskip\labelsep\ignorespaces\normalfont}
{\hfill\openbox\par\addvspace{6pt}}
\def\b{\backslash}
\def\k{\langle k \rangle}
\def\tr{\textrm}
\def\la{\langle}
\def\ra{\rangle}
\def\alp{\alpha^{(r-k+1)}}

\def\A{\mathcal{A}}
\def\B{\mathcal{B}}
\def\C{\mathcal{C}}
\def\D{\mathcal{D}}
\def\E{\mathcal{E}}
\def\F{\mathcal{F}}
\def\mc{\mathcal}
\newtheorem{theorem}{Theorem}

\newtheorem{lemma}{Lemma}                    
\newtheorem{construct}{Construction}

\begin{document}
\title{The number of $k$-intersections of an intersecting family of $r$-sets}
\author{John Talbot
\\Merton College
\\University of Oxford
\\E-mail: talbot@maths.ox.ac.uk}
\date{\today}
\maketitle
\begin{abstract}  
The Erd\H os-Ko-Rado theorem tells us how large an intersecting family of $r$-sets from an $n$-set can be, while results due to Lov\'asz and Tuza give bounds on the number of singletons that can occur as pairwise intersections of sets from such a family.

We consider a natural common generalization of these problems. Given an intersecting family of $r$-sets from an $n$-set and $1\leq k \leq r$, how many $k$-sets can occur as pairwise intersections of sets from the family? For $k=r$ and $k=1$ this reduces to the problems described above. We answer this question exactly for all values of $k$ and $r$, when $n$ is sufficiently large. We also characterize the extremal families.
\end{abstract}
\section{Introduction}
A family of sets is \emph{intersecting} if any two sets from the family meet. Let $[n]^{(r)}$ denote the collection of all $r$-sets from $[n]=\{1,2,\dots ,n\}$, $[n]^{(\leq r)}$ denote the collection of all sets of size at most $r$ from $[n]$ and $2^{[n]}$ denote the power-set of $[n]$. For an intersecting family $\A\subseteq 2^{[n]}$ and $1\leq k \leq n$ define the \emph{intersection structure} of $\A$ by $$I(\A)=\{A\cap B:A,B\in\A\}$$ and the collection of \emph{$k$-intersections} of $\A$ by $$\A\langle k \rangle=\{A\in I(\A):|A|=k\}.$$
Note that although $\A\la 1\ra$ is a collection of singleton sets it will often be convenient to treat it as a set of points, for example $\{a,b,c\}$ instead of $\{\{a\},\{b\},\{c\}\}$.

The primary result concerning intersecting families of sets is the celebrated Erd\H os-Ko-Rado theorem. This tells us exactly how large an intersecting family $\A\subseteq [n]^{(r)}$ can be. In the notation introduced above it gives a bound on $|\A\langle r\rangle|$.
\begin{theorem}[Erd\H os, Ko and Rado \cite{EKR}]
\label{EKRthm}
Let $n\geq 2r$ and $\mathcal{A} \subset [n]^{(r)}$ be intersecting. Then $|\mathcal{A}|=|\A\langle r \rangle| \leq \binom{n-1}{r-1} =|\mathcal{A}_{1}|$, where $\mathcal{A}_{1}=\{A\in [n]^{(r)}:1\in A\}$. 
\end{theorem} 
This theorem has been generalized in many different ways. Amongst the most significant results in this area are the Hilton-Milner theorem \cite{HILMIL}, the Ray-Chaudhuri-Wilson theorem \cite{RCW}, the Hajnal-Rothschild theorem \cite{HAJROT} and the Complete Intersection theorem \cite{AHLKHA}. 

All of these results share the common aim of giving bounds on the size of a family of sets satisfying certain intersection properties. In this paper we take a different approach. Rather than giving bounds on the size of an intersecting family itself, we instead consider the intersection structure of such a family. In particular we consider the question: how many sets of a given size can occur as the intersection of two sets from an intersecting family? In other words, given an intersecting family $\A\subseteq 2^{[n]}$ and $1\leq k\leq n$ how large can $|\A\la k \ra|$ be?

If $\A\subseteq 2^{[n]}$ is intersecting and $1\leq k \leq n$ then clearly $\A\langle k \rangle\subseteq [n]^{(k)}$. In fact we may have $\A\langle k \rangle= [n]^{(k)}$. (Consider, for $n$ odd, the collection of all sets in $2^{[n]}$ containing more than $n/2$ points.) However if we restrict ourselves to uniform families this question becomes more interesting.

We define for $1\leq k \leq r \leq n$ \[
\beta(n,r,k)=\max\{|\A\langle k \rangle|:\A\subseteq [n]^{(r)}\textrm{ is intersecting}\}.\]
Theorem \ref{EKRthm} may now be rephrased as: ``if $n\geq 2r$ then $\beta(n,r,r)=\binom{n-1}{r-1}$''.

It is perhaps not immediately obvious that $\beta(n,r,1)$ is bounded above for fixed $r$, irrespective of the value of $n$. However the following result, due to Lov\'asz \cite{LOV}, shows that this is indeed the case.
\begin{theorem}[Lov\'asz \cite{LOV}]
\label{LOVthm}
If $r\geq 1$ and $$ \alpha^{(r)}=\max\{\beta(n,r,1):n \geq r\}$$ then $\alpha^{(r)}$ is well-defined and satisfies \begin{equation*}\binom{2r-3}{r-1} +2r-2\leq\alpha^{(r)} \leq (2r-1)\binom{2r-3}{r-1}.\end{equation*}
\end{theorem}
Theorem \ref{LOVthm} was subsequently improved by Tuza \cite{TUZA1} who gave the following bounds for $\alpha^{(r)}$.
\begin{theorem}[Tuza \cite{TUZA1}]\label{TUZAthm}
If $r\geq 4$ then $$2\binom{2r-4}{r-2}+2r-4\leq \alpha^{(r)} \leq \binom{2r-1}{r-1}+\binom{2r-4}{r-1}.$$
\end{theorem}
The lower bound in Theorem \ref{TUZAthm} comes from the following construction. Take $[2r-4]^{(r-2)}=\{A_1,B_1,\ldots,A_m,B_m\}$, with $m=\frac{1}{2}\binom{2r-4}{r-2}$ and $A_i\dot{\cup} B_i=[2r-4]$. For each $1\leq i \leq m$ introduce four new vertices: $a_i,b_i,c_i$ and $d_i$. Then our intersecting family $\A$ consists of the sets $A_i\cup\{a_i,b_i\},A_i\cup\{c_i,d_i\},B_i\cup\{a_i,c_i\}$ and $B_i\cup\{b_i,d_i\}$ for $i=1.\ldots,m$. Clearly $\A\la 1 \ra$ contains $[2r-4]$ together with $\{a_i,b_i,c_i,d_i:1\leq i \leq m\}$.

Although Theorem \ref{TUZAthm} gives bounds on $\alpha^{(r)}$ that are sharp up to a multiplicative constant factor, $\alpha^{(r)}$ is only known exactly for $r\leq 4$, with $\alpha^{(1)}=1$, $\alpha^{(2)}=3$, $\alpha^{(3)}=7$ and $\alpha^{(4)}=16$.
\section{Main Result}
Our main result, Theorem \ref{mainthm}, is an exact determination of $\beta(n,r,k)$ for all $1\leq k \leq r$ and sufficiently large $n$. This result says that in order to maximize the number of $k$-intersections of an intersecting family of $r$-sets we should first take an intersecting family of $(r-k+1)$-sets whose pairwise intersections realize as many singletons as possible and then extend this to an intersecting family of $r$-sets by taking all $r$-sets containing a member of this family.

Since the value of $\alpha^{(r)}$ is only known for $r\leq 4$ our answer is necessarily given in terms of $\alpha^{(r-k+1)}$.

In order to describe our main result we need to introduce a construction, which we will show is the essentially unique extremal family. However before we can give this construction we require a lemma.

We say that an intersecting family $\A\subseteq [n]^{(r)}$ is \emph{maximal} iff any set in $\A\b [n]^{(r)}$ is disjoint from at least one set in $\A$.
\begin{lemma}\label{maxlem} If $n\geq 2r$ and $\A \subseteq [n]^{(r)}$ is intersecting and maximal then for every $A,B\in\A$ we have
\[
A\cap B\cap \A\la 1\ra\neq \emptyset.\]
\end{lemma}
\begin{proof}
Suppose $A,B\in\A$ and $A\cap B=\{c_1,\ldots,c_k\}=C$ is disjoint from $\A\la 1\ra$. Let $D=\{d_1,\ldots d_k\}\subseteq [n]\b (A\cup B)$. Now, since each $c_i\not\in\A\la 1\ra$ and using the maximality of $\A$, if  we replace each $c_i$ in $A$ by the corresponding $d_i$ then the resulting set will still belong to $\A$. Hence $(A\b C) \cup D \in \A$, but this is disjoint from $B$, a contradiction. 
\end{proof}
Throughout the remainder of this section we will assume that $n$ is a large positive integer without explicitly determining exactly how large it must be. We will return to this question at the end of the next section. At this point we simply remark that in order to give our construction we will require at least $n\geq \alp$. 
\begin{construct}\label{example}
Let $1\leq k\leq r <n$, with $n$ large. Let $\mc{B}\subseteq [n]^{(r-k+1)}$ be an intersecting family satisfying $\B\langle 1 \rangle=[\alp]$. Define a new intersecting family 
\[
\A=\{A\in [n]^{(r)}:\exists B\in\mc{B}\textrm{ such that } B \subseteq A\}.\]
We claim that
$$\A\langle k \rangle=\{C\in [n]^{(k)}:C\cap [\alp]\neq \emptyset\}.$$
\end{construct} 
\begin{proof}[Proof of Claim]
We show first that any two sets from $\A$ contain a common point in $[\alp]$. By adding sets from $[n]^{(r-k+1)}$ to $\B$ we may form a maximal intersecting family $\C$ such that $\B\subseteq \C\subseteq [n]^{(r-k+1)}$. Since $\B\la 1 \ra \subseteq \C \la 1\ra$ and $|\B\la 1 \ra|=\alp$ we have $\C\la 1\ra=\B\la 1 \ra=[\alp]$. 

If $B_1,B_2\in\B\subseteq \C$ then by Lemma \ref{maxlem} we have $B_1\cap B_2\cap [\alp]\neq \emptyset$. So if $A_1,A_2\in\A$ then there exist $B_1,B_2\in \B$ with $B_1\subseteq A_1$ and $B_2\subseteq A_2$. Hence $A_1\cap A_2\cap [\alp]\neq \emptyset$. Thus 
\[
\A\k \subseteq \{C\in [n]^{(k)}:C\cap [\alp]\neq \emptyset\}.\]

We now prove the other inclusion.

Let $C\in [n]^{(k)}$ satisfy $C\cap [\alp]\neq \emptyset$. If $x \in C \cap [\alp]$ then there exist $B_1,B_2\in\mc{B}$ such that $B_1\cap B_2=\{x\}$. Let $D_1=B_1\cup C$ and $D_2=B_2\cup C$. If $|D_1|=d_1$ and $|D_2|=d_2$ let $E_1$ and $E_2$ be disjoint sets of sizes $r-d_1$ and $r-d_2$ respectively in $[n]\b (D_1\cup D_2)$. Then $A_1=D_1\dot{\cup} E_1$ and $A_2=D_2\dot{\cup} E_2$ both belong to $\A$ and $A_1\cap A_2=C$. So $C \in \A\langle k \rangle$ as claimed. 
\end{proof}

\begin{figure}[!ht]
\begin{center}
\psfrag{A}[][]{\Large$\A$}
\psfrag{B}[][]{\Large$\mathcal{B}=[3]^{(2)}$}
\psfrag{C}[][]{$\cdots$}
\includegraphics[width=250pt]{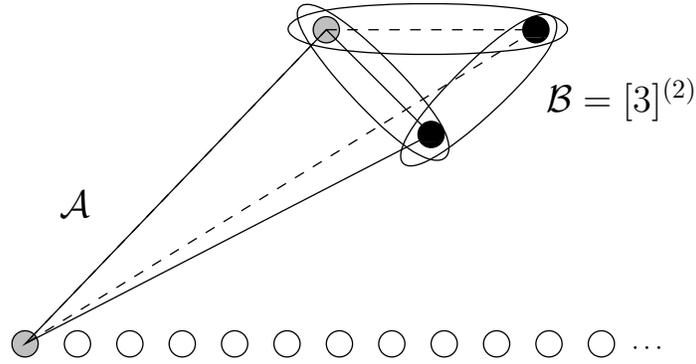}
\caption{Detail of an extremal family for $k=2$ and $r=3$}
\label{intfig}
\end{center}
\end{figure}

 Figure \ref{intfig} attempts to show how part of this family looks for $k=2$ and $r=3$.

Construction \ref{example} shows that for $n$ sufficiently large
\[
\beta(n,r,k)\geq |\{A\in [n]^{(k)}: A\cap [\alp]\neq \emptyset\}|.\]
Our main result says that in fact this is best possible, and that any other family achieving this bound may be constructed in a similar fashion.
\begin{theorem}
\label{mainthm}
If $1\leq k \leq r <n$, with $n$ large,  then \[
\beta(n,r,k)= |\{A\in [n]^{(k)}: A\cap [\alp]\neq \emptyset\}|=\alp\binom{n}{k-1}+O(n^{k-2}).\]Moreover if $\A\subseteq [n]^{(r)}$ is intersecting and maximal then either 
\[
|\A\k|\leq (\alp-1)\binom{n}{k-1}+O(n^{k-2}),\]
or there is a family $\C\subseteq [n]^{(r-k+1)}$ satisfying $|\C\la 1 \ra|=\alp$ such that 
\[
\B=\{B\in [n]^{(r)}:\exists C\in \C\tr{ such that }C\subseteq B\}\]
is contained in $\A$ and \[
\A\la l \ra =\{A\in [n]^{(l)}:A\cap \C\la 1 \ra\neq \emptyset\},\]
for $k\leq l \leq r$.
\end{theorem}
\section{Proof}
We prove Theorem \ref{mainthm} using the following two lemmas. The first says that either $|\A\k|$ is small or $\A$ contains a family similar to that given in Construction \ref{example}. The second lemma then allows us to show that we cannot extend such a family so as to increase $|\A\k|$.
\begin{lemma}\label{conlem} If $\A\subseteq [n]^{(r)}$ is intersecting and maximal then either\[
|\A\k|\leq (\alp-1)\binom{n}{k-1}+O(n^{k-2}),\]
or there exists a set $D\in [n]^{(k-1)}$ such
\[
\C=\{A\b D: D\subset A \in \A\}\]
satisfies $|\C\la 1\ra|=\alp$ and
\[
\B=\{B\in [n]^{(r)}:\exists C\in \C\tr{ such that }C\subseteq B\}\]
is contained in $\A$.
\end{lemma}
\begin{lemma}\label{extlem} If $1\leq k\leq l\leq r \leq n$ and $A,B\in [n]^{(r)}$, $\C\subseteq [n]^{(\leq r-k+1)}$ satisfy:
\begin{itemize}
\item[(i)] $A\cap B=D\in [n]^{(l)}$,
\item[(ii)] $\E=\C\cup\{A,B\}$ is intersecting,
\item[(iii)] $|\C\la 1 \ra|=\alp$,
\end{itemize}
then $D\cap\C\la 1 \ra\neq \emptyset$.
\end{lemma}
\begin{proof}[Proof of Theorem \ref{mainthm}]
Construction \ref{example} implies that for $n$ large
\[
\beta(n,r,k)\geq |\{A\in [n]^{(k)}: A\cap [\alp]\neq \emptyset\}|.\]
The fact that this is also an upper bound for $\beta(n,r,k)$ will follow if we prove the remainder of Theorem \ref{mainthm}. 

So let $\A\subseteq [n]^{(r)}$ be intersecting and maximal. Lemma \ref{conlem} then implies that either \[
|\A\k|\leq (\alp-1)\binom{n}{k-1}+O(n^{k-2}),\]
or there exists a set $B\in [n]^{(k-1)}$ such
\[
\C=\{A\b B: B\subset A \in \A\}\]
satisfies $|\C\la 1\ra|=\alp$ and
\[
\B=\{B\in [n]^{(r)}:\exists C\in \C\tr{ such that }C\subseteq B\}\]
is contained in $\A$. We may suppose that the later holds. Hence, as in Construction \ref{example}, we have 
\[
\B\k=\{A\in [n]^{(k)}:A\cap \C\la 1\ra\neq \emptyset\}\subseteq \A\k.\]
Indeed, for any $1\leq i \leq r$ we have \[
\B\la i \ra =\{A\in [n]^{(i)}:A\cap \C\la 1\ra\neq\emptyset\}\subseteq \A\la i \ra.\]
In order to complete the proof of Theorem \ref{mainthm} it is sufficient to show that if $k \leq l \leq r$ and $D\in\A\la l\ra$ then $D\cap \C\la 1\ra \neq \emptyset$.

 If $D\in\A\la l \ra$ then there exist $A,B\in \A$ such that $A\cap B=D\in [n]^{(l)}$. Also $\E=\C\cup\{A,B\}$ is intersecting, since if $F\in\A$ is disjoint from $C\in\C$ then there exists $G\in\B\subseteq \A$ containing $C$ such that $F\cap G=\emptyset$, a contradiction. So $A,B,\C,D$ and $\E$ satisfy the conditions of Lemma \ref{extlem} and hence $D\cap \C\la 1 \ra \neq \emptyset$ as required.
\end{proof}
We now turn to the proofs of Lemmas \ref{conlem} and \ref{extlem}. 
\begin{proof}[Proof of Lemma \ref{conlem}]
Let $\A\subseteq [n]^{(r)}$ be intersecting and maximal. If $a=|\A\la 1\ra|<\alp$ then Lemma \ref{maxlem} implies that
\begin{eqnarray*}
|\A\k|&\leq&|\{A\in [n]^{(k)}:A\cap \A\la 1 \ra\neq \emptyset\}|\\
& \leq &(\alp-1)\binom{n}{k-1}+O(n^{k-2}).\end{eqnarray*}
So we may suppose that $a \geq \alp$.

We partition $\A\k$ as $\A\k=\A_1\dot{\cup}\A_2$, where
\[
\A_1=\{A\in\A\k:\tr{if $D\subset A$ and $|D|=k-1$ then $D\cap \A\la 1\ra\neq \emptyset$}\}\]
and $\A_2=\A\k \b \A_2$.

If $A\in\A_1$ then $|A\cap \A\la 1\ra|\geq 2$, since otherwise there exists $D\subset A$ satisfying $|D|=k-1$ and $D\cap \A\la 1\ra= \emptyset$. Moreover since $a\leq \alpha^{(r)}=O(1)$ we have
\[
|\A_1|\leq \sum_{i=2}^{k}\binom{a}{i}\binom{n-a}{k-i}=O(n^{k-2}).\]
For $D\in [n]^{(k-1)}$ define
\[
\A_D=\{A\b D: D\subset A \in \A\}.\] 
Then
\begin{equation}\label{sum}
|\A_2|\quad\leq \sum_{D\in [n]^{(k-1)},\, D\cap \A\la 1\ra=\emptyset}|\A_D\la 1\ra|.\end{equation}
If $D\in [n]^{(k-1)}\b\A \la k-1\ra$ then $\A_D\subset [n]^{(r-k+1)}$ is intersecting. Also, by Lemma \ref{maxlem}, if $D\in \A\la k-1 \ra$ then $D\cap \A\la 1\ra\neq \emptyset$. Hence each term in the sum (\ref{sum}) is bounded above by $\alp$. So either there exists a set $D\in [n]^{(k-1)}$ such that $D\cap \A\la 1\ra=\emptyset$ and $|\A_D\la 1\ra|=\alp$ or
\[
|\A_2|\leq (\alp-1)\binom{n-a}{k-1}.\]
In the latter case we have
\[
|\A\k|=|\A_1|+|\A_2|\leq (\alp-1)\binom{n}{k-1} + O(n^{k-2}),\]
so we may suppose that the former holds. 

Let $D\in [n]^{(k-1)}$ satisfy $D\cap \A\la 1\ra=\emptyset$ and $|\A_D\la 1\ra|=\alp$. Define 
\[
\C=\A_D=\{A\b D:D\subset A \in \A\}.\]
Since $|\C\la 1 \ra|=\alp$ the proof will be complete if we can show that
\[
\B=\{B\in [n]^{(r)}:\exists C\in \C\tr{ such that }C\subseteq B\}\subseteq \A.\] If $B\in\B$ then there exists $C\in\C$ such that $C\subseteq B$. By definition of $\C$ we have $C\cup D\in\A$ and $D\cap\A\la 1\ra=\emptyset$. Now $C$ meets every set in $\A$, since if $A\in \A$ and $A\cap C=\emptyset$ then $A\cap(C\cup D)\cap \A\la 1\ra=\emptyset$, contradicting Lemma \ref{maxlem}. Hence $B$ meets every set in $\A$ and so by maximality $B\in\A$.
\end{proof}
\begin{proof}[Proof of Lemma \ref{extlem}] We use induction on $k$, for $1 \leq k \leq r$. If $k=1$ then $\E=\C\cup \{A,B\}$ is an intersecting family in $[n]^{(r)}$. Let $\E'$ be a maximal intersecting family in $[n]^{(r)}$ containing $\E$. Then $\C\la 1 \ra\subseteq \E'\la 1\ra$ and $|\C\la 1\ra|=\alpha^{(r)}$ imply that $\E'\la 1\ra=\C\la 1\ra$. So by Lemma \ref{maxlem} $A\cap B \cap \C\la 1\ra \neq \emptyset$ as required. Hence the result holds for $k=1$.

Let $2\leq k \leq r$ and assume the result holds for $k-1$. Suppose, for a contradiction, that $D\cap \C\la 1 \ra=\emptyset$. Since $2\leq k\leq l$, there exist $a,b\in D\b \C\la 1 \ra$ with $a\neq b$. We now replace $a$ and $b$ by $v^*\not\in [n]$. Define $A^*=A\b\{a,b\} \cup \{v^*\}$, $B^*=B\b\{a,b\} \cup \{v^*\}$ and $D^*=D\b\{a,b\} \cup \{v^*\}$. Also let $\C^*,\E^*$ be the families produced from $\C,\E$ respectively by replacing each occurrence of $a$ or $b$ in every set in these families by $v^*$. Note that $\C^*$ and $\E^*$ are intersecting and $\C^*$ is a family of sets of size at most $r-k+1$.

If $\C^*\la 1\ra=\C\la  1 \ra$ then, since $1\leq k-1\leq l-1\leq r-1$, we may apply our inductive hypothesis for $k-1$ to $A^*,B^*,\C^*,D^*$ and $\E^*$, with $l,r$ replaced by $l-1,r-1$ respectively. This implies that $D^*\cap \C^*\la 1 \ra\neq \emptyset$. Then $\C^*\la 1\ra=\C\la 1 \ra$ and $v^*\not\in \C\la 1\ra$ imply that $D\cap \C\la 1 \ra\neq \emptyset$ as required. So we may suppose that $\C^*\la 1 \ra\neq \C\la 1 \ra$.

Since $\C^*$ is an intersecting family of sets of size at most $r-k+1$, $|\C\la 1 \ra|=\alp$ and $\C^*\la 1 \ra\neq \C\la 1 \ra$, there must exist $c\in\C\la 1 \ra$ such that $c\not\in \C^*\la 1 \ra$. This means that in replacing $a,b$ by $v^*$ we have ``lost'' the intersection $\{c\}$. This can only happen if for all $F,G\in \C$
\begin{equation}\label{star}
F\cap G=\{c\}\implies (a\in F\tr{ and }b\in G) \tr{ or }(b\in F\tr{ and }a\in G).\end{equation}
Take such a pair $F,G\in \C$ (they exist since $c\in\C\la 1 \ra$). Without loss of generality we may suppose that $a\in F$ and $b\in G$. Consider the set $F\b\{c\}\cup \{b\}$. We claim that this set meets every set in $\C$. If not then there exists $H\in\C$ such that $H\cap(F\b\{c\}\cup\{b\})=\emptyset$. Thus $F\cap H=\{c\}$ and $b\not\in H$ so by (\ref{star}) we must have $a\in H$. Hence $a\in F\cap H$, contradicting the fact that $F\cap H=\{c\}$. 

So $\F=\C\cup \{F\b\{c\}\cup \{b\}\}$ is an intersecting family of sets of size at most $(r-k+1)$. However $(F\b\{c\}\cup\{b\})\cap G=\{b\}$ and $b\not\in\C\la 1 \ra\subset \F\la 1 \ra$ imply that $|\F\la 1 \ra|\geq \alp+1$, contradicting the definition of $\alp$. The result then follows by induction. 
\end{proof}
\section{Remarks}
Can the characterization of the extremal families in Theorem \ref{mainthm} be strengthened to $\A=\B$ rather than $\B\subseteq \A$ and $\A\la l \ra=\B\la l\ra$ , for $k\leq l \leq r$ (using the notation of Theorem \ref{mainthm})? For $k=1$ or $k=r$ this is trivially true, while for $k=r-1$ this can still be verified easily. However in general this is false.

For example, if $k=2$ and $r=4$, let
\[
\D=\{123,145,246,356,167,257\}\subset [7]^{(3)}\]
and
\[
\E=\{E\in [n]^{(4)}:\exists D\in\D\tr{ such that }D\subset E\}.\]
Then $|\D\la 1\ra|=\alp=7$ and the family
\[
\F=\E\cup\{1258,3478\}\]
satisfies \[
\F\la 2 \ra=\{F\in [n]^{(2)}: F\cap [7]\neq \emptyset\}.\]
However it is easy to check that any set in $[n]^{(4)}\b \F$ that meets every set in $\F$ must belong to $[8]^{(4)}$. Hence if $\A$ is a maximal intersecting family in $[n]^{(4)}$ containing $\F$ then the family $\B$ given by Theorem \ref{mainthm} must be $\D$. So in this case we not only have $\A\neq \B$ but also $\A\la k-1\ra\neq \B\la k-1\ra$, since $\B\la 1\ra=[7]\neq [8]=\A\la 1\ra$.

Another obvious question to ask is when can we actually evaluate the expression given in Theorem \ref{mainthm} for particular values of $k,r$ and $n$? As we remarked earlier, $\alpha^{(r)}$ is only known for $r\leq 4$. So we can evaluate $\alp$ and hence \[
|\{A\in [n]^{(k)}:A\cap [\alp]\neq \emptyset\}|\]
for $r-k+1\leq 4$. For other values of $k$ and $r$ we can use Theorem \ref{TUZAthm} to give bounds.

Finally we turn to the question of how large $n$ must be for the value of $\beta(n,r,k)$ to be determined by Theorem \ref{mainthm}. In order to use Construction \ref{example} to give a lower bound for $\beta(n,r,k)$ we require $n\geq \alp$ and so Theorem \ref{TUZAthm} implies we need\[ n\geq 2\binom{2r-2k-2}{r-k-1}+2r-2k-2.\]

Conversely, examining the proof of Lemma \ref{conlem} our argument requires
\[
\sum_{i=1}^k\binom{\alp}{i}\binom{n-\alp}{k-i}\geq (\alp -1)\binom{n-a}{k-1}+\sum_{i=2}^{k}\binom{a}{i}\binom{n-a}{k-i},\]
where $a=|\A\la 1\ra|$, for the upper bound to be valid. So, since $\alp\leq a \leq \alpha^{(r)}$, a rough calculation shows that Theorem \ref{mainthm} determines $\beta(n,r,k)$ for $n\geq k(\alpha^{(r)}+1)^2$. Theorem \ref{TUZAthm} then implies that $\beta(n,r,k)$ is determined for \[n\geq k\binom{2r}{r}^2.\]

Clearly this could be improved, indeed it is plausible that Construction \ref{example} is best possible whenever it exists. However, given our lack of knowledge of the true value of $\alpha^{(r)}$ it seems an extremely difficult problem to determine $\beta(n,r,k)$ exactly for small values of $n$.


\begin{thebibliography}{}
\bibitem{AHLKHA} R. Ahlswede and L.H. Khachatrian, The complete intersection theorem for systems of finite sets, \emph{European J. Combin.} \textbf{18} (1997), 125-136.
\bibitem{BOL} B. Bollob\'as, On generalized graphs, \emph{Acta Math. Acad. Sci. Hungar.} \textbf{16} (1965) 447-452.
\bibitem{ERD} P. Erd\H os, A problem on independent $r$-tuples, \emph{Ann. Univ. Sci. Budapest} \textbf{8} (1965), 93-95.
\bibitem{EKR} P. Erd\H os, C. Ko and R. Rado, Intersection theorems for systems of finite sets, \emph{Quart. J. Math. Oxford (2)} \textbf{12} (1961), 313-320.
\bibitem{HAJROT} A. Hajnal and B. Rothschild, A generalization of the Erd\H os-Ko-Rado theorem on finite sets, \emph{J. Combin. Theory Ser. A} \textbf{15} (1973), 359-362.
\bibitem{HILMIL} A.J.W. Hilton and E.C. Milner, Some intersection theorems for systems of finite sets, \emph{Quart. J. Math. Oxford (2)} \textbf{18} (1967), 369-384.
\bibitem{LOV} L. Lov\'asz, ``Combinatorial Problems and Exercises'', North Holland, Amsterdam, New York, Oxford, 1979.
\bibitem{RCW} D.R. Ray-Chaudhuri and R.M. Wilson, On $t$-designs, \emph{Osaka J. Math.} \textbf{12} (1975), 737-744.
\bibitem{TUZA1} Zs. Tuza, Critical hypergraphs and intersecting set-pair systems, \emph{J. Combin. Theory Ser. B}, \textbf{39} (1985), 134-145.
\bibitem{TUZA2} Zs. Tuza, Applications of the set-pair method in extremal hypergraphs, in: P. Frankl et al. (eds.), ``Extremal Problems for Finite Sets'', Bolyai Society Mathematical Studies, Vol 3 (1994), 479-514,  Janos Bolyai Math. Society.
\end{thebibliography}
\end{document}